\documentclass[12pt]{amsart}

\usepackage{amsfonts, amsmath, amssymb}
\usepackage{graphicx, graphics}

\newtheorem{prop}{Proposition}[section]
\newtheorem{theo}[prop]{Theorem}

\newtheorem{conj}[prop]{Conjecture}

\theoremstyle{definition}
\newtheorem{dfn}[prop]{Definition}
\newtheorem{exam}[prop]{Example}

\theoremstyle{remark}
\newtheorem{rem}[prop]{Remark}

\def\<{\langle}
\def\>{\rangle}

\def\Si{{\Sigma}}

\def\P{{\mathbb P}}
\def\R{{\mathbb R}}
\def\T{{\mathbb T}}
\def\C{{\mathbb C}}
\def\Z{{\mathbb Z}}
\def\Q{{\mathbb Q}}

\def\A{{\mathcal A}}

\def\N{{\widetilde{N}}}
\def\M{{\widetilde{M}}}
\def\dd{{\widetilde{d}}}
\def\nn{{\mathbf n}}
\def\D{{\widetilde{\Delta}}}

\begin{document}

\title{Mixed toric residues and Calabi-Yau complete intersections}

\author{Victor~V.~Batyrev}
\address{Mathematisches Institut \\ Universit\"at T\"ubingen \\
Auf der Morgenstelle 10 \\ T\"ubingen D-72076, Germany} 
\email{victor.batyrev@uni-tuebingen.de} 

\author{Evgeny~N.~Materov}
\address{Mathematisches Institut \\ Universit\"at T\"ubingen \\
Auf der Morgenstelle 10 \\ T\"ubingen D-72076, Germany} 
\email{evgeny.materov@uni-tuebingen.de}



\begin{abstract}
Using the Cayley trick, we define the notions of mixed toric residues and 
mixed Hessians associated with $r$ Laurent polynomials $f_1, \ldots, f_r$.
We conjecture that the values of mixed toric residues on the  mixed Hessians 
are determined by mixed volumes of the Newton polytopes of $f_1, 
\ldots, f_r$. Using mixed toric residues, we generalize our 
Toric Residue Mirror Conjecture to the case of 
Calabi-Yau complete intersections in Gorenstein toric Fano varieties 
obtained from nef-partitions of reflexive polytopes. 
\end{abstract}

\maketitle

\tableofcontents

\newpage

\section{Introduction}

This paper  is the continuation of our previous work  \cite{Batyrev-Materov1} 
where 
we proposed a  toric mirror symmetry test using toric residues. The idea 
of this test  has appeared in the paper of Morrison-Plesser 
\cite{MP} who  observed that the coefficients of some 
power series expansions of unnormalized Yukawa couplings for mirrors of 
 Calabi-Yau 
hypersurfaces in toric varieties $\P$ can be interpreted as generating  
functions for  intersection numbers of divisors on some sequences of 
toric varieties $\P_{\beta}$ parameterized by lattice points $\beta$ in the 
Mori cone $K_{\rm eff}(\P)$ of $\P$. Due to results of 
Mavlyutov \cite{Mavlyutov2}, 
it is known that the unnormalized Yukawa couplings 
can be computed using the toric residues introduced by Cox \cite{Cox}. 
 In our  paper \cite{Batyrev-Materov1}, we formulated 
a  general mathematical conjecture, so called {\it 
Toric Residue Mirror Conjecture}, which describes 
some  power series expansions of the toric residues in terms of intersection 
numbers of divisors on a sequence of simplicial toric varieties 
$\P_{\beta}$ (we call them Morrison-Plesser moduli spaces).
This  conjecture includes  all examples of mirror symmetry for Calabi-Yau
hypersurfaces in Gorenstein toric varieties associated with reflexive 
polytopes. Since the toric mirror symmetry construction exists  also 
for Calabi-Yau complete intersection in Gorenstein toric Fano varieties 
\cite{Borisov,Batyrev-Borisov1}, 
it is natural to try to extend our conjecture to this more general 
situation.
\medskip
     
The case of Calabi-Yau complete intersections of $r$  
hypersurfaces $$f_1(t) = \cdots = f_r(t)= 0, \;\; r>1,$$ 
defined by  Laurent polynomials $f_1(t),\ldots,f_r(t)
\in\C[t_1^{\pm 1},\ldots,t_d^{\pm 1}]$ in $d$-dimensional toric 
varieties $\P$ was not considered by Morrison and Plesser in \cite{MP}. 
We remark that in this case one {\it does not} get a connection to the 
``quantum cohomology ring'' \cite{BatyrevQuantum} as in  
the hypersurface case. This difference is explained by the consideration of 
a non-reflexive $(d+r-1)$-dimensional polytope  
$\D$, so called Cayley polytope, and its   secondary polytope 
${\rm Sec}(\D)$. The Cayley polytope  $\D$ 
appears from the Cayley trick which introduces $r$   additional $r$ variables 
$t_{d+1}, \ldots, t_{d+r}$ and a new polynomial 
$F(t):= \sum_{j=1}^r t_{d+j}f_j(t)$. 
We consider the usual toric residue ${\rm Res}_F$  
associated  with $F$ and define  
the  {\it $k$-mixed toric residue}  ${\rm Res}^k_F$ corresponding to  
a positive integral 
solution $k = (k_1, \ldots, k_r)$ 
of the equation $k_1 + \cdots + k_r = d+r$ as a $k$-th homogeneous 
component of  ${\rm Res}_F$.   
We expect that the  $k$-mixed toric residues are  similar to the usual 
toric residues. In particular, we introduce the notion of  
the {\it $k$-mixed Hessian $H^k_F$} of Laurent  polynomials  
$ f_1, 
\ldots, f_r $ and conjecture that the 
value of  ${\rm Res}^k_F$ on $H_F^k$ is exactly 
the mixed volume 
\[ 
 V(\underbrace{\Delta_1, \ldots, \Delta_1}_{k_1-1}, 
\ldots,\underbrace{\Delta_r, \ldots, \Delta_r}_{k_r-1}), \]
where $\Delta_1, \ldots, \Delta_r$ are Newton polytopes of $f_1, \ldots, 
f_r$. 
\medskip

Our generalization of the Toric Residue Mirror Conjecture for Calabi-Yau
complete intersections uses the notions of the  nef-partition $\Delta =
\Delta_1 + \cdots + \Delta_r$ of $d$-dimensional 
reflexive polytope $\Delta$ \cite{Borisov}. 
In this situation, one obtains  
a dual nef-partition $\nabla = \nabla_1 + \cdots + \nabla_r$ and two more 
reflexive polytopes: 
\[ \nabla^* = conv\{ \Delta_1, \ldots, \Delta_r \}, \;\;  
\Delta^* = conv\{ \nabla_1, \ldots, \nabla_r \} .\] 
It is important that special coherent triangulations of $\nabla^*$ define 
coherent triangulations of the Cayley polytope 
$\D := \Delta_1 * \cdots * \Delta_r$. Therefore  
the choice of  such a triangulation
${\mathcal T}$ of  $\nabla^*$ determines a vertex $v_{\mathcal T}$ of 
the secondary polytope ${\rm Sec}(\D)$ and a partial projective simplicial 
crepant desingularization $\P : = \P_{\Sigma({\mathcal T})}$ of the Gorenstein
toric variety $\P_{\nabla^*}$. So 
one obtains a sequence of simplicial toric varieties $\P_\beta$ associated 
with the lattice points $\beta$ in the Mori cone $K_{\rm eff}(\P)$ of $\P$.
 We conjecture that 
the generating function of intersection numbers 
\[ I_P(a) = \sum_{\beta \in K_{\rm eff}(\P)}  I(P, \beta)\, a^\beta \] 
coincides with the power series expansion of the $k$-mixed toric residue
\[ R_P(a) = {\rm Res}_F^k(P(a,t)) \] 
at the vertex $v_{\mathcal T} \in  {\rm Sec}(\D)$.  The precise formulation 
of this conjecture is given in Section 4. 

In Sections 5, 6 we check our conjecture for nef-partitions corresponding 
to Calabi-Yau complete 
intersections in weighted projective spaces 
$\P(w_1,\ldots,w_n)$
 and in product of projective spaces $\P^{d_1}\times\cdots\times\P^{d_p}$. 
The final section is devoted to 
applications of the Toric Residue Mirror Conjecture to the 
computation of Yukawa couplings 
for Calabi-Yau complete intersections. 

\medskip

{\it Acknowledgments.} This work was supported by DFG, Forschungsschwerpunkt
``Globale Methoden in der komplexen Geometrie''. E.~Materov was supported in
part by RFBR Grant 00-15-96140. We are also appreciate to hospitality 
and support 
of the Isaac Newton Institute in Cambridge.

\bigskip

\section{Toric residues} 

In this section we remind necessary well-known facts about toric residues 
(see \cite{Cox,CDS,Batyrev-Materov1}).

Let $\M$ and $\N = {\rm Hom}(\M,\Z)$ be two free abelian groups of rank 
$\dd$ dual to each other. We denote by 
\[
  \<*,*\>\,:\, \M\times \N\rightarrow \Z
\]
the natural bilinear pairing, and by $\M_\R$ (resp. $\N_\R$) the real 
scalar extension of $\M$ (resp. $\N$).

\begin{dfn}[\cite{Batyrev-Borisov2}] 
A $\dd$-dimensional rational polyhedral cone $C$ ($\dd > 0$) in $\M_\R$ 
is called {\it Gorenstein} if it is strongly convex (i.e.,  
$C + (-C) = \{0\}$), there exists an element 
$n_C\in \N$ such that $\<x,n_C\> > 0$ for any nonzero $x\in C$, and all
vertices of the $(\dd - 1)$-dimensional convex polytope
\[
  \Delta(C) = \{x\in C\,:\,\<x,n_C\> = 1\}
\]
belong to $\M$. The polytope $\Delta(C)$ is called the 
{\it supporting polytope of $C$}. 
For any $m\in C\cap \M$, we define the {\it degree of $m$} as 
\[
  \deg m = \<m,n_C\>.
\]
\end{dfn}

\begin{dfn}
\label{pol-ring}
Let $\D = \Delta(C)$ be the  supporting polytope for a  Gorenstein cone 
$C\subset \M_\R$. We denote by  $S_\D$ 
the  semigroup $\C$-algebra of the monoid of lattice points $C\cap \M$. 
In order to transform  the additive semigroup operation in  $C\cap \M$ into a 
multiplicative form in  $S_\D$, we write $t^m$ for the element 
in $S_\D$ corresponding to $m \in C$. One can consider $S_\D$ as a 
graded $\C$-algebra:
\[
  S_\D = \bigoplus_{l = 0}^\infty S_\D^{(l)},
\]
where  
the $l$-th homogeneous component $S_\D^{(l)}$ has a  
$\C$-basis consisting of all $t^m$ such that  $m\in C\cap \M$ and 
$\deg m = l$.
We define also the homogeneous ideal 
\[
  I_\D = \bigoplus_{l = 0}^\infty I_\D^{(l)}
\]
in $S_\D$ whose  $\C$-basis consists of 
all $t^m$ such that $m$ is a lattice point in the interior  of 
$C$. 
\end{dfn}

\begin{dfn} 
An element 
\[ g:= \sum_{m \in \D \cap \M} a_m t^m\in S_\D^{(1)}, 
\;\; a_m \in \C\]
is called $\D$-regular if for some $\Z$-basis $n_1, \ldots, n_{\dd}$ of $\N$
the elements 
\[ g_i:= \sum_{m \in \D \cap \M} a_m \< m, n_i \> t^m,\;\; i=1, \ldots, 
\dd \]
form a regular sequence in $S_\D$. We define the 
matrix $G:= (g_{ij})_{1 \leq i, j \leq \dd} $, where 
\[   g_{ij}:= \sum_{m \in \D \cap \M} a_m \< m, n_i \> \< m, n_j \>  
t^m,\;\; i,j=1, \ldots, 
\dd. \] 
The element
\[ H_g := \det G \]
is called {\it Hessian of $g$}.
\end{dfn} 

\begin{rem} 
a) The definition of $\D$-regularity does not depend on the choice 
of $\Z$-basis $n_1, \ldots, n_{\dd}$ of $\N$. It gives the same 
non-degeneracy condition as the one proposed by Mavlyutov in 
\cite{Mavlyutov1}.   
In many applications the lattice vector $n_C$ will be included  in 
$\{ n_1, \ldots, n_{\dd}\}$. \\
b) If $\D \cap \M = 
\{ m_1, \ldots, m_{\mu} \}$, then by 
\cite[Proposition~1.2]{CDS}, one has 
\[ H_g = \sum_{ 1 \leq i_1 < \ldots < i_{\dd} \leq \mu}  
\left( \det(m_{i_1}, \ldots, m_{i_{\dd}}) \right)^2 t^{m_{i_1} +  \cdots + 
m_{i_\dd}} .\] 
In particular, $H_g$ is independent on the choice of the 
$\Z$-basis $n_1, \ldots, n_{\dd}$ and $H_g \in I_\D^{(\dd)}$.\\ 
c) The graded $\C$-algebra  $S_\D$ is Cohen-Macaulay and 
$I_\D$ is its dualizing module. If $g$ is $\D$-regular in 
$S_\D$, then 
\[
  S_g := S_\D/\<g_1,\ldots,g_{\dd}\>S_\D,
\]
is a  graded finite-dimensional ring and 
\[
  I_g := I_\D/\<g_1,\ldots,g_{\dd}\>I_\D
\]
is a  graded  $S_g$-module together with 
a non-degenerate pairing 
\[
  S_g^{(l)} \times I_g^{(\dd - l)}\rightarrow I_g^{(\dd)}\simeq \C, \quad
  l = 0,\ldots, \dd - 1.
\]
induced by the $S_g$-module structure.
\label{rem-res}
\end{rem} 

\begin{dfn} 
By the {\it toric residue} corresponding to a $\D$-regular element 
$g \in S^{(1)}_\D$ we mean the $\C$-linear mapping 
\[ {\rm Res}_g \; : \; I_\D^{(\dd)} \to \C \]
which is uniquely determined  by two conditions: 

(i) ${\rm Res}_g(h) = 0$ for any $h \in \< g_1, \ldots, g_{\dd} \> I_\D$; 

(ii) ${\rm Res}_g(H_g) = {\rm Vol}(\D)$, 
where  ${\rm Vol}(\D)$ denotes the volume of the $(\dd-1)$-dimensional 
polytope $\D$ multiplied by $(\dd -1)!$. 
\label{t-residue}
\end{dfn} 

Let $\P_\D:= {\rm Proj}\, 
S_\D$ be the $(\dd-1)$-dimensional toric variety associated 
with the polytope $\D$ and ${\mathcal O}_{\P_\D}(1)$ the corresponding ample
sheaf on $\P_\D$. Then one has the canonical isomorphisms of graded rings
\[ S_\D \cong \bigoplus_{l \geq 0} H^0(\P_\D,{\mathcal O}_{\P_\D}(l)) \]
and graded modules
\[  I_\D \cong \bigoplus_{l \geq 0} H^0(\P_\D,{\omega}_{\P_\D}(l)), \] 
where ${\omega}_{\P_\D}$ is the dualizing sheaf on  $\P_\D$. In particular, 
we obtain a canonical isomorphism 
\[  I_\D^{(\dd)} \cong  H^0(\P_\D,{\omega}_{\P_\D}(\dd)). \]

The following statement is a simple reformulation of Theorem 2.9(i) 
in \cite{Batyrev-Materov1}: 

\begin{prop} 
 Let $n_1, \ldots, n_{\dd}$ be  a $\Z$-basis of $\N$ such that 
$n_1 = n_C$. Denote by $m_1, \ldots, m_{\dd}$ the dual $\Z$-basis of 
$\M$. For any  elements $h \in  I_\D^{(\dd)}$ and 
 $g \in   S_\D^{(1)}$, we  define  a rational differential 
$(\dd -1)$-form on $\P_\D$: 
\[ \Omega(h,g):= \frac{ h}{g_1 \cdots g_{\dd}} \frac{dt^{m_2}}{t^{m_2}} 
\wedge \cdots \wedge  \frac{dt^{m_{\dd}}}{t^{m_{\dd}}}. \]
If $g$ is   $\D$-regular, then 
\[ {\rm Res}_g(h) = \sum_{\xi \in V_g}  {\rm res}_{\xi} 
  \left( \Omega(h,g) \right),
\]
where $V_g = \{ \xi \in \P_\D \, : \, g_2(\xi) = 
\cdots = g_\dd(\xi) = 0 \}$  is 
the set of common zeros of $g_2, \ldots, g_{\dd}$ and 
${\rm res}_{\xi}(\Omega(h,g))$ is the local Grothendieck residue of the
form $\Omega(h,g)$ at the point $\xi \in V_g$. 

In particular, if all the common roots of $g_2,\ldots,g_{\dd}$  are simple 
and contained in the open dense $(\dd -1)$-dimensional torus $\T 
\subset \P_\D$,
then 
\[  {\rm Res}_g(h)
   =  
  \sum_{\xi \in V_g}
  \frac{ p(\xi)}{g_1(\xi) {H}_g^1(\xi)}, 
\]
where  ${H}_g^1$ is the determinant of the matrix $G^1:= 
\left( g_{ij}) \right)_{ 2 \leq i,j \leq \dd}$. 
\label{sum_res}
\end{prop}

\begin{dfn}[\cite{Batyrev-Borisov1}]
A Gorenstein cone $C$ is called {\it reflexive} if the 
dual cone 
\[
  \check{C} = \{y\in\N_\R\,:\,\<x,y\>\ge 0\quad \forall x\in C\}
\]
is also Gorenstein, i.e., there exists 
$m_{\check{C}}\in \M$ such that $\<m_{\check{C}},y\> > 0$ 
for all $y\in \check{C}\setminus\{0\}$, and all vertices of the 
supporting polytope 
\[
  \Delta(\check{C}) = \{y\in \check{C}\,:\,\<m_{\check{C}},y\> = 1\}
\]
belong to $\N$. We will call the integer 
$r = \<m_{\check{C}},n_C\>$ the {\it index of} $C$ (or $\check{C}$).
A $(\dd - 1)$-dimensional lattice polytope 
$\D$ is called {\it reflexive} if it is a supporting polytope 
of some $\dd$-dimensional reflexive Gorenstein cone $C$ of index $1$.
Moreover, the supporting polytope $\D^*$ of 
the dual cone $\check{C}$ is also reflexive polytope which is called 
{\it dual} (or {\it polar}) to $\D$.
\end{dfn}

If $C$ is a reflexive Gorenstein cone of index $r$, then $I_\D$ is 
a principal ideal generated by the element $t^{m_{\check{C}}}$ of degree $r$. 
So one obtains the 
canonical isomorphism $I_\D^{(l)} \cong S_\D^{(l-r)}$. In particular,
there exists the toric residue mapping
\[ {\rm Res}_g\: : \; S_\D^{(\dd-r)} \to \C\] 
which   
is uniquely determined  by the  conditions: 

(i) ${\rm Res}_g(h) = 0$ for any $h \in \< g_1, \ldots, g_{\dd} \> S_\D$; 

(ii) ${\rm Res}_g(H_g') = {\rm Vol}(\D)$,
where $H_g = t^{m_{\check{C}}}H_g'$. 

\bigskip

\section{The Cayley trick and mixed toric residues}

Let $M$ be a free abelian group of rank $d$, $M_\R:= M \otimes \R$, and 
$\Delta \subset M_\R$ a convex $d$-dimensional polytope with vertices 
in $M$. We assume that there exist $r$ convex polytopes $\Delta_1, 
\ldots, \Delta_r$ with  vertices in $M$ such that $\Delta$ 
can be written as the  Minkowski sum $\Delta = 
\Delta_1 + \cdots + \Delta_r$ (here we do not require  that all polytopes 
$\Delta_1, \ldots, \Delta_r$ have maximal dimension $d$).     

\begin{dfn} 
We set $\M :=  M \oplus \Z^r$, $\dd := d + r$ and define the $\dd$-dimensional
Gorenstein 
cone $C=C(\Delta_1, \ldots, \Delta_r)$ in $\M_\R :=  M_\R \oplus \R^r$ 
as follows
\[
  C := \{(\lambda_1 x_1 + \cdots + \lambda_r x_r, \lambda_1,\ldots,\lambda_r)
  \in \M_\R\,:\,\lambda_i \ge 0,\, x_i \in \Delta_i,\,
  i = 1, \ldots, r\}. 
\]
The $(d+r -1)$-dimensional polytope  
$\Delta_1* \cdots *\Delta_r$  
 defined as the intersection of the cone $C$ with the affine hyperplane 
$\sum_{i = 1}^r \lambda_i = 1$
\begin{eqnarray*}
\Delta_1* \cdots *\Delta_r:= 
 \{ ( 
  \lambda_1 x_1 + \cdots + \lambda_r x_r, \lambda_1, \ldots, \lambda_r)
 \,:\,\lambda_i \ge 0,\, 
  \sum_{i = 1}^r \lambda_i = 1,\, x_i \in \Delta_i\},
\end{eqnarray*}
will be called the  
{\it Cayley polytope associated with the Minkowski sum decomposition 
$\Delta = \Delta_1 + \cdots + \Delta_r$}.  
It is clear that all vertices of  $\Delta_1* \cdots *\Delta_r$ are contained 
in $\M$ and 
\begin{eqnarray*}
\Delta_1* \cdots *\Delta_r   =
 conv((\Delta_1 \times \{b_1\})\cup\cdots\cup(\Delta_r \times \{b_r\})
), 
\end{eqnarray*}
where $\{b_1,\ldots,b_r\}$ is the standard basis of $\Z^r$. For 
fixed polytopes $\Delta_1, \ldots, \Delta_r$ we 
denote $\Delta_1* \cdots *\Delta_r$ simply by $\D$. 
\end{dfn} 

\begin{dfn} 
Define $S_{\D} : =\C[C\cap \M]$ to be the semigroup algebra of the monoid 
$C\cap \M$ over the complex numbers. The algebra   $S_{\D}$ has a natural 
$\Z_{\geq 0}^r$-grading defined by the last $r$ coordinates of lattice points 
in $\M$. By choosing an isomorphism $M \cong \Z^d$, we can identify 
$S_{\D}$ with a $\Z_{\geq 0}^r$-graded monomial subalgebra in 
$$\C[t_1^{\pm 1},
\ldots, t_d^{\pm 1}, t_{d+1}, \ldots, t_{d+r} ], $$
where the $\Z_{\geq 0}^r$-grading is considered with respect to the last $r$ 
variables 
$t_{d+1}, \ldots, t_{d+r}$. We denote by $I_\D$ the $\Z_{\geq 0}^r$-graded 
monomial ideal in $S_\D$ generated by all lattice points in the interior 
of $C$. For any $k=(k_1, \ldots, k_r) \in \Z_{\geq 0}^r$, we denote by
$S^k_\D$  (resp. by $I_\D^k$) the $k$-homogeneous component of $S_\D$ 
(resp. of  $I_\D$). We will use also the total  $\Z_{\geq 0}$-grading on 
$S_\D$ and  $I_\D$. For any nonnegative integer $l$, we denote  the 
corresponding $l$-homogeneous components of $S_\D$ and  $I_\D$  by $S_\D^{(l)}$ and 
 $I_\D^{(l)}$ respectively. So one has:
\[ S_\D^{(l)} = \bigoplus_{|k| = l} S^k_\D, \;\; \;\; 
 I_\D^{(l)} = \bigoplus_{|k| = l} I^k_\D, \]
where  $|k|:= k_1 + \cdots + k_r$.
\end{dfn} 

Let $f_1(t), \ldots, f_r(t)$ be Laurent polynomials in 
$\C[t_1^{\pm 1}, \ldots, t_d^{\pm 1} ]$ such that $\Delta_i$ is 
the Newton polytope of $f_i$ $(1 \leq i \leq r)$. We set  
\[ F(t) := t_{d+1}f_1(t) + \cdots + t_{d+r} f_r(t). \]
It is easy to see  that $\D = \Delta_1* \cdots *\Delta_r$ is the Newton
polytope of $F$. Moreover, using the decomposition 
\[  S_\D^{(1)} = \bigoplus_{|k| = 1} S^k_\D =  \bigoplus_{i = 1}^r 
 S^{b_i}_\D ,\] 
we see  that every Laurent polynomial $G$ 
in $\C[t_1^{\pm 1}, 
\ldots, t_d^{\pm 1}, t_{d+1}, \ldots, t_{d+r} ] $
with the Newton polytope $\D$  can be obtained from the sequence 
of arbitrary Laurent polynomials 
$g_1, \ldots, g_r \in \C[t_1^{\pm 1}, \ldots, t_d^{\pm 1} ]$ 
by the formula $G =  t_{d+1}g_1 + \cdots + t_{d+r} g_r$, where 
$\Delta_i$ is the Newton polytope of $g_i$ $(1 \leq i \leq r)$.
The above correspondence  $\{ f_1, \ldots, f_r \}  \mapsto F$ is 
usually called the {\it Cayley trick}. We call $F =   t_{d+1}f_1 + 
\cdots + t_{d+r} f_r$ the {\it Cayley polynomial associated with 
$f_1, \ldots, f_r$}. 

\begin{dfn} 
Let $\Delta_1, \ldots, \Delta_r \subset M_\R$ be convex polytopes 
with  vertices in $M$ such $\Delta = \Delta_1 + \cdots + \Delta_r$ has 
dimension $d$. We say that $r$  Laurent polynomials 
\[ f_i(t) = \sum_{m \in \Delta_i\cap M} a_m^{(i)}t^m, \;\; 
i = 1, \ldots,  r \]
form a {\it $\D$-regular sequence} if the corresponding Cayley polynomial
$F$ is  $\D$-regular, i.e., the polynomials
\[ F_i := t_i\partial/\partial t_i F, \;\; i =1, \ldots, d+ r \]
form a regular sequence in $S_\D$. 
\end{dfn}

\begin{dfn} 
Let $f_1(t), \ldots, f_r(t)$ be Laurent polynomials with Newton
polytopes $\Delta_1, \ldots, \Delta_r$ as above, $F =   t_{d+1}f_1 + 
\cdots + t_{d+r} f_r$ the corresponding Cayley polynomial, and 
\[  H_F := 
  \det\left(t_i\frac{\partial F_j}
  {\partial t_i}\right)_{1 \le i,j \le d + r} = 
  \det\left(\left(t_i\frac{\partial}{\partial t_i}\right)
            \left(t_j\frac{\partial}{\partial t_j}\right)F
  \right)_{1 \le i,j \le d + r} \in I^{(d+r)}_\D \subset 
S^{(d+r)}_\D  \]
the Hessian of $F$. For any $k = (k_1, \ldots, k_r)$ with 
$|k| = d+r$ we define $H_F^k \in I^k_\D$ to be the $k$-homogeneous component 
of $H_F$. The polynomial  $H^k_F$ will be called {\it $k$-mixed Hessian 
of $f_1, \ldots, f_r$. } 
\end{dfn} 

\begin{rem} 
Since the last $r$ rows of the matrix 
\[ \left(\left(t_i\frac{\partial}{\partial t_i}\right)
            \left(t_j\frac{\partial}{\partial t_j}\right)F
  \right)_{1 \le i,j \le d + r} \]
are divisible respectively by $t_{d+1}, \ldots, t_{d+r}$, the Hessian 
$H_F$ is divisible by  the monomial  $t_{d+1}\cdots t_{d+r}$.  
Therefore  $H_F^k =0$ if one of the coordinates $k_i$  of $k = 
(k_1, \ldots, k_r)$ 
is  zero. 
In particular, one has 
\[ H_F = \sum_{\begin{subarray}{c}
k \in \Z_{>0}^r \\  |k| = d+r
\end{subarray}} H^k_F.
\]
\end{rem} 

Let $k = (k_1, \ldots, k_r) \in \Z_{>0}^r$ be a solution of the 
linear Diophantine equation
\[ |k|=  k_1 + \cdots + k_r = d + r. \]
For any  $r$ subsets $S_i \subset \Delta_i \cap M$  
such that $|S_i| = k_i$ $( 1\leq i \leq r)$ we define the nonnegative 
integer $\nu(S_1, \ldots, S_r)$ as follows: choose an element $s_i$ 
in each $S_i$ $(1\leq i\leq r)$, define $S$ to be  
the $d\times d$-matrix whose rows  are 
all possible nonzero vectors $s- s_i$ where $s \in S_i$, 
$1\leq i\leq r$, and set $\nu(S_1, \ldots, S_r) := (\det S)^2$. It is easy 
to see that up to sign $\det S$ does not depend on the choice 
of elements $s_i \in S_i$ and  therefore  $\nu(S_1, \ldots, S_r)$ is well 
defined. 

\begin{prop} 
Let $k = (k_1, \ldots, k_r) \in \Z_{> 0}^r$ be a positive integral  solution of the 
linear Diophantine equation
\[ |k|=  k_1 + \cdots + k_r = d + r. \]
Then the mixed Hessian can be computed by the following formula 
\[ H^k_F = t_{d+1}^{k_1} \cdots t_{d+r}^{k_r} 
\sum_{(S_1, \ldots, S_r)}  
\nu(S_1, \ldots, S_r) \prod_{i =1}^r \prod_{s_i \in S_i} a_{s_i}^{(i)} 
t^{s_i}, \]
where the sum runs over all $r$-tuples $(S_1, \ldots, S_r)$ of 
 subsets $S_i \subset \Delta_i \cap M$  
such that $|S_i| = k_i$ $( 1\leq i \leq r)$. 
\label{k-mixed}
\end{prop} 

\begin{proof} The  formula for $H_F^k$ follows immediately from 
the formula in \ref{rem-res}(b) 
applied to the Cayley polytope $\D$. 
\end{proof} 

\begin{dfn} 
Let $k = (k_1, \ldots, k_r) \in \Z_{> 0}^r$ be a positive integral 
 solution of the equation 
\[ |k|=  k_1 + \cdots + k_r = d + r. \]
Consider the toric residue  
\[ {\rm Res}_F \; : \; I_\D^{(d + r)}\to \C\] defined
as a $\C$-linear map which vanishes on $\< F_1, \ldots, F_{d+r} \> I_\D$ 
and sends 
$H_F$ to ${\rm Vol}(\D)={\rm Vol}(\Delta_1* \cdots *\Delta_r)$. 
The restriction 
 ${\rm Res}_F^k$  of 
${\rm Res}_F$ to the $k$-th homogeneous component  $I_\D^{k}$:
\[ {\rm Res}_F^k \; : \;  I^k_\D \to \C \]
 will be called  the 
{\it $k$-mixed toric 
residue associated with $f_1, \ldots, f_r$.}
\end{dfn}

Since $H_F^k$ is an element of $I^k_\D$, it is natural to ask about 
the value of ${\rm Res}_F^k(H_F^k)$. 

\begin{conj}
Let $k = (k_1, \ldots, k_r) \in \Z_{> 0}^r$ be a positive integral 
 solution of $ k_1 + \cdots + k_r = d + r$. We set 
$\overline{k} = (\overline{k_1}, \ldots, \overline{k_r}) := 
(k_1 -1, \ldots, k_r -1)$. 
Then 
\[ {\rm Res}_F^k(H_F^k) =
 V(\underbrace{\Delta_1, \ldots, \Delta_1}_{\overline{k_1}}, 
\ldots,\underbrace{\Delta_r, \ldots, \Delta_r}_{\overline{k_r}}), \]
where $V(\Theta_1,\ldots, \Theta_d)$ denotes the mixed volume of convex 
polytopes $\Theta_1, \ldots, \Theta_d$ multiplied by $(d + r - 1)!$. 
\label{normaliz}
\end{conj}

Our conjecture agrees with a  result of Danilov and Khovanskii:  

\begin{prop} \cite[\S 6]{Danilov-Khovanskii} 
The normalized volume of the Cayley polytope $\D=
\Delta_1* \cdots *\Delta_r$ can be computed by the following formula:
\[ {\rm Vol}(\Delta_1* \cdots *\Delta_r)= \sum_{|\overline{k}| = d} 
 V(\underbrace{\Delta_1, \ldots, \Delta_1}_{\overline{k_1}}, 
\ldots,\underbrace{\Delta_r, \ldots, \Delta_r}_{\overline{k_r}}). \]
\end{prop} 

\begin{rem} 
Let  $r =d$ and $k = (d+1, 1, \ldots, 1)$. It follows from 
\ref{k-mixed} and  \ref{rem-res}(b)  that  
\[ H_F^k = H_{t_{d+1}f_1}(t_{d+2}f_2) \cdots (t_{2d}f_d), \]
where 
$$H_{t_{d+1}f_1}= \det\left(\left(t_i\frac{\partial}{\partial t_i}\right)
            \left(t_j\frac{\partial}{\partial t_j}\right)t_{d+1}f_1
  \right)_{1 \le i,j \le d + 1}.  $$ 
On the other hand, we have
\[  V(\underbrace{\Delta_1, \ldots, \Delta_1}_{\overline{k_1}}, 
\ldots,\underbrace{\Delta_r, \ldots, \Delta_r}_{\overline{k_r}})
=  V(\underbrace{\Delta_1, \ldots, \Delta_1}_d) =
{\rm Vol}(\Delta_1). \]
Therefore, Conjecture \ref{normaliz} can be considered as a 
``generalization'' of  \ref{t-residue}(ii). 
\end{rem} 
\medskip

It is easy to show that 
the cone $C$ from 3.1 is a reflexive Gorenstein cone of index $r$ 
if and only if $\Delta = \Delta_1 + \cdots + \Delta_r$ is a reflexive
polytope.  
In this situation,   we have  
\[ I_\D =  t_{d+1} \cdots  t_{d+r}S_\D. \]
Therefore one has canonical isomorphisms:
\[ I_\D^k \cong  S^{\overline{k}}_\D, \;\;\forall\;  k \in \Z^r_{>0},\]
where  the monomial basis in  $S^{\overline{k}}_\D$ can be identified 
with the set of all lattice points in $\overline{k_1}\Delta_1 + 
\cdots + \overline{k_r}\Delta_r$.  The ${\overline{k}}$-homogeneous 
component of the corresponding toric residue map 
\[ {\rm Res}_F^{\overline{k}} \; : \;  S^{\overline{k}}_\D \to \C. \]
will be also called the {\it  ${\overline{k}}$-mixed toric residue}.
 
\bigskip

\section{Toric Residue Mirror Conjecture}

Let $M$ and $N = {\rm Hom}(M,\Z)$ be the dual to each other abelian groups  
of rank $d$, $M_\R$ and $N_\R$ their  
$\R$-scalar extensions and $\Delta \subset M_\R$ 
a reflexive polytope with the unique interior lattice point $0 \in M$. 
Denote by $\P_{\Delta}$ the  Gorenstein toric Fano variety associated with  
$\Delta$. Let $D_1,\ldots,D_s$ be the toric divisors on $\P_\Delta$ 
corresponding to the codimension-1 faces $\Theta_1, \ldots, \Theta_s$ and 
$e_1, \ldots, e_s$ the vertices of the dual reflexive polytope 
$\Delta^* \subset N_{\R}$ such that 
\[ \Delta = \{ x \in M_\R \; : \; \< x, e_j \> \geq -1, \;\; 
j = 1, \ldots, s \} , \] 
\[ \Theta_j = \Delta \cap \{  x \in M_\R \; : \; \< x, e_j \> = -1 \}, 
\;\;j \in \{ 1, \ldots, s \}.   \] 

\begin{dfn}
A  Minkowski sum  
$\Delta = \Delta_1 + \cdots + \Delta_r$ is called a {\it nef-partition} 
of the reflexive polytope $\Delta$ 
if all vertices of $\Delta_1, \ldots, \Delta_r$ belong to $M$, 
and 
\[ 
\min_{x \in \Delta_i } \< x, e_j \> \in \{ 0, -1 \}, \;\; 
\forall 1 \leq i \leq r, \; \forall 1 \leq j \leq s.\]
\end{dfn}

Since $\min_{x \in \Delta }  \< x, e_j \> = -1$ for all 
$j \in \{ 1, \ldots, s \}$, the equality 
$\min_{x \in \Delta_i } \< x, e_j \> =-1$ 
holds exactly for one index $i \in \{ 1, \ldots, r \}$ if 
we fix a vertex $e_j \in \Delta^*$. Therefore, we 
can split the set of vertices 
$ \{ e_1, \ldots, e_s \} \subset \Delta^*$ into a disjoint union of subsets 
$B_1, \ldots, B_r$ where 
\[ B_i: = \{ e_j \in \{1, \ldots, s \} \; : 
\;   \min_{x \in \Delta_i } \< x, e_j \> =-1\} . \]
Now we can define $r$ nef Cartier divisors
\[  E_i := \sum_{j \, : \, e_j \in B_i} D_j, \;\; i = 1,\ldots,r. \]
Therefore, a   nef-partition $\Delta = \Delta_1 + \cdots + \Delta_r$ of
polytopes 
induces a partition of the anti-canonical divisor 
$-K_{\P_\Delta} = D_1 + \cdots + D_n$ of $\P_{\Delta}$ 
into a sum of $r$ nef Cartier divisors: 
\[ -K_{\P_\Delta} = E_1 + \cdots + E_r. \] 

Now it is easy to see that the above definition of the nef-partition is 
equivalent to the definition given in \cite{Borisov}. 

\begin{dfn} 
If $\Delta = \Delta_1 + \cdots + \Delta_r$ is a nef-partition, then 
for any $i = 1,\ldots,r$ we denote 
\[ \nabla_i := \{ y \in N_\R\; : \; \< x, y \> \geq - \delta_{ij}, \;\;
x \in \Delta_j, \; j = 1, \ldots,  r \}. \]
The lattice polytopes $\nabla_1, \ldots, \nabla_r$ define another
nef-partition $\nabla := \nabla_1 + \cdots + \nabla_r$ of the  reflexive 
polytope $\nabla \subset N_\R$ which is called the  
{\it dual nef-partition}. 
\end{dfn}

The lattice polytopes $\nabla_1, \ldots, \nabla_r$ can be also defined as
\[
  \nabla_j := conv(\{0\}\cup B_j)\subset M_\R,\quad j = 1,\ldots,r.
\]
Moreover, one has two dual reflexive polytopes 
\[ \Delta^* = conv(\nabla_1 \cup \cdots \cup \nabla_r) \subset 
N_\R \] 
\[ \nabla^* = conv(\Delta_1 \cup \cdots \cup \Delta_r) \subset 
M_\R. \]  

A nef-partition $\Delta = \Delta_1 + \cdots + \Delta_r$  defines a family of 
$(d - r)$-dimensional Calabi-Yau complete 
intersections defined by vanishing of 
$r$ Laurent polynomials $f_1, \ldots, f_r$ with Newton polytopes 
$\Delta_1, \ldots, \Delta_r$. 
According to \cite{Borisov}, the dual nef-partition 
$\nabla = \nabla_1 + \cdots + \nabla_r$ defines 
the  mirror dual family of Calabi-Yau complete intersections.
\medskip

Define  $\A_j$ to be a subset in 
$\Delta_j \cap M$ containing all vertices of $\Delta_j$  and set  
$A_j := \A_j \setminus \{ 0 \}$,  $ j =1, \ldots, r$. It is easy to 
see that $A_i \cap A_j = \emptyset$ for all $ i \neq j$.  
We set  
$$A_1 \cup \cdots \cup A_r := 
\{ v_1, \ldots, v_n\},$$ 
and define $a_1, \ldots, a_n\in\C$ 
to be the coefficients
of the Laurent polynomials $f_1, \ldots, f_r$ 
\[ f_j(t) := 1 - \sum_{i\; : \; v_i \in A_j} a_i t^{v_j} ,  
\;\; j =1, \ldots, r.\]
Let $A := \{ 0 \} \cup A_1\cup \cdots \cup A_r$ and 
$\D = \Delta_1 * \cdots * \Delta_r$ 
be  the Cayley 
polytope. Denote by 
$\pi$ the injective mapping  
\[   A_1\cup \cdots \cup A_r \to \D \cap \M  \] 
which sends  a nonzero lattice point $m \in A_j$ to $(m,b_j)$ 
$(1 \leq j \leq r)$ and define
\[ \widetilde{A}: = \pi ( A_1\cup \cdots \cup A_r) \cup  
\{ (0,b_1), \ldots, 
(0,b_r) \}. \]

We hold notations 
from \cite[\S 4]{Batyrev-Materov1}.

\begin{dfn} 
 Choose  a coherent triangulation 
${\mathcal T} = \{ \tau_1, \ldots, \tau_p \}$ of the reflexive polytope 
$\nabla^* = conv(\Delta_1 \cup \cdots \cup \Delta_r)$ associated with $A$ 
such that $0$ is a vertex of all its $d$-dimensional simplices 
$ \tau_1, \ldots, \tau_p$ and all elements of $A$ appear as vertices of 
simplices in ${\mathcal T}$.  Define  a coherent triangulation 
$\widetilde{\mathcal T} =  \{ \widetilde{\tau_1}, \ldots, 
 \widetilde{\tau_p} \}$ of $\D = \Delta_1 * \cdots * \Delta_r$ 
associated with $ \widetilde{A}$ as follows: a 
$(d+r-1)$-dimensional simplex $ \widetilde{\tau_i} \in 
\widetilde{\mathcal T}$ 
is the convex hull of $\pi$-images of all nonzero vertices 
of $\tau$ and $\{ (0,b_1), \ldots, 
(0,b_r) \}$. We call $\widetilde{\mathcal T}$ the {\bf induced triangulation} 
of $\D$. 
\end{dfn}

 Let $\P:= \P_{\Sigma({\mathcal T})}$ be the $d$-dimensional
simplicial toric variety defined by the fan $\Sigma({\mathcal T}) \subset 
M_\R$ ($\P$ is a partial crepant desingularization of the Gorenstein
toric Fano variety $\P_{\nabla}$) and denote by $\P_{\beta}$ the
 {\it Morrison-Plesser 
moduli space} \cite[Definition~3.3]{Batyrev-Materov1} 
corresponding to a lattice point 
$$\beta = (\beta_1, \ldots, \beta_n) \in R(\Sigma) = \{ (x_1, \ldots, x_n) 
\in \Z^n\; : \; x_1 v_1 + \cdots + x_n v_n = 0 \} $$   
in the Mori cone $K_{\rm eff}(\P)$. 
One has a canonical surjective homomorphism 
$$\psi_{\beta}\, : \, H^2(\P, \Q) \to H^2 (\P_{\beta}, \Q).$$

\begin{dfn} 
By abuse of notations, let us denote by $[D_j] \in  H^2 (\P_{\beta}, \Q)$ 
$( 1 \leq j \leq n)$ the image of $[D_j] \in  H^2(\P, \Q)$ under $\psi_\beta$. 
Using the multiplication in the cohomology ring 
$H^* (\P_{\beta}, \Q)$, we define the intersection product
 \[ \Phi_\beta
 := 
  [E_1]^{(E_1,\beta)}\cdots [E_r]^{(E_r,\beta)}
  \prod_{i:(D_i,\beta) < 0} [D_i]^{-(D_i,\beta) - 1}
\]
considered as a cohomology class in $H^{2({\rm dim} \, \P_{\beta} - d)} 
( \P_{\beta}, \Q)$ and call $\Phi_\beta$ the {\bf Morrison-Plesser class}
corresponding to the nef-partition $\Delta = \Delta_1 + \cdots + \Delta_r$.  
\end{dfn} 

\begin{dfn} 
Let $k = (k_1, \ldots, k_r) \in \Z_{> 0}^r$ be a positive integral 
 solution of 
\[ |k|=  k_1 + \cdots + k_r = d + r. \] 
A polynomial $P(x_1, \ldots, x_n) \in \Q[x_1, \ldots, x_n]$ is called
{\bf $\overline{k}$-homogeneous} if it is homogeneous of degree 
$\overline{k_i} = k_i -1$ with respect 
to every group of $|A_i|$  variables $x_j$ $( v_j \in A_i)$ $( 1 \leq i 
\leq r)$. 
\end{dfn} 

Now we are able to formulate a generalized Toric Residue Mirror Conjecture:

\begin{conj} 
\label{TRMC} Let $\Delta = \Delta_1 + \cdots + \Delta_r$ and 
 $\nabla = \nabla_1 + \cdots + \nabla_r$ be two arbitrary 
dual nef-partitions. 
 Choose any  coherent triangulation 
${\mathcal T} =\{\tau_1, \ldots, \tau_p\} $ of $\nabla^*$ associated with $A$
such that $0$ is a vertex of all the simplices 
$\tau_1, \ldots, \tau_p$ as above. 
Then for any $\overline{k}$-homogeneous polynomial 
$P(x_1,\ldots,x_n)\in\Q[x_1,\ldots,x_n]$
of degree $d$ the Laurent expansion of the ${\overline{k}}$-mixed
 toric residue 
\[
  R_P(a) := (-1)^d\,{\rm Res}_F^{\overline{k}}
(t_{d+1}^{\overline{k_1}} \cdots t_{d+r}^{\overline{k_r}}\, 
P(a_1 t^{v_1},\ldots,a_n t^{v_n}))
\]
at the vertex $v_{\widetilde{\mathcal T}} \in {\rm Sec}(\D)$ 
corresponding to the induced 
triangulation $\widetilde{{\mathcal T}}=  \{ \widetilde{\tau_1}, \ldots, 
 \widetilde{\tau_p} \}$ coincides with 
the generating function of 
intersection numbers 
\[
  I_P(a) := \sum_{\beta\in K_{\rm eff}(\P)} I(P, \beta)\, a^\beta, 
\]
where the sum runs over all integral points $\beta = (\beta_1,\ldots,
\beta_n)$ of the
Mori cone $K_{\rm eff}(\P)$, $a^\beta := a_1^{\beta_1}\cdots a_n^{\beta_n}$, 
\[
  I(P, \beta) = \int_{\P_\beta}P([D_1],\ldots,[D_n])\Phi_\beta = 
  \< P([D_1],\ldots,[D_n])\Phi_\beta \>_\beta, 
\]
and $\Phi_\beta \in H^{2(\dim \P_\beta -d)}(\P_{\beta}, \Q)$ is 
the Morrison-Plesser class of $\P_{\beta}$. We assume $I(P,\beta)$ 
to be zero if $\P_\beta$ is empty. 
\end{conj}
\bigskip

\section{Complete intersections in weighted projective spaces}

Let $\P = \P(w_1,\ldots,w_n)$ be a $d$-dimensional weighted
projective space, $n = d + 1$. The fan $\Si$ of $\P(w_1,\ldots,w_n)$
is determined by $n$ vectors $v_1,\ldots,v_n\in M\simeq \Z^d$ which 
generate $M$ and  
satisfy the relation
\[
  w_1 v_1 + \cdots + w_n v_n = 0.
\]
If we assume that  ${\rm gcd}(w_1,\ldots,w_n) = 1$ and
\[
  w_i | (w_1 + \cdots + w_n), \quad i = 1,\ldots,n,
\]
then $\P$ is a Gorenstein toric Fano variety with the anticanonical 
divisor $-K_\P = D_1 + \cdots + D_n$, where $D_i$ is the 
toric divisor corresponding to the vector $v_i$.
These divisors are related modulo rational equivalence as
\[
  \frac{[D_1]}{w_1} = \cdots = \frac{[D_n]}{w_n} =: [D_0].
\]
Consider a decomposition $\{v_1,\ldots,v_n\}$
into a disjoint union of $r$ nonempty 
subsets $A_1,\ldots,A_r$ and 
define  the divisors  $E_i:= \sum_{j\, : \, v_j \in A_i} D_j$ 
 on $\P$ such that 
$[E_i] = d_i [D_0]$, where $d_j = \sum_{i \in A_j} w_i$, $j = 1,\ldots,r$.
Note that the integers $d_i$ satisfy 

$d_1 + \cdots + d_r = w_1 + \cdots + w_n$. Let $\Delta_i := 
conv( \{ 0\} \cup A_i)$ $(1 \leq i \leq r)$. The polytopes 
$\Delta_1, \ldots, \Delta_r$ 
define a nef-partition $\Delta:= \Delta_1 + \cdots + \Delta_r$ 
if and only if 
\[  w_i | d_j, \quad i = 1,\ldots,n,\; j =1, \ldots, r.
\]

The following  result generalizes \cite[Theorem 7.3]{Batyrev-Materov1}:

\begin{theo} Let $P\in \Q[x_1,\ldots,x_n]$ be a homogeneous polynomial of
degree $d$. Then the generating function of intersection numbers on the
Morrison-Plesser moduli spaces has the form
\[
  I_P(y) = \nu \cdot P(w_1,\ldots,w_n) \sum_{b\ge 0} \mu^b y^b =
  \frac{\nu \cdot P(w_1,\ldots,w_n)}{1 - \mu y},
\]
where 
\[
  \nu := \frac{1}{w_1 \cdots w_n},\quad
  \mu := \frac{d_1^{d_1} \cdots d_r^{d_r}}{w_1^{w_1} \cdots w_n^{w_n}}, \quad
 y := a_1^{w_1}\cdots a_n^{w_n}.
\]
\end{theo}

\begin{proof} The lattice points $\beta$ in the Mori cone of $\P$ correspond
to the linear relations $b w_1 v_1 + \cdots + bw_n v_n = 0$, $b \in \Z_{\geq 0}$.
Therefore we set  $y := a_1^{w_1}\cdots a_n^{w_n}$. 

The  Morrison-Plesser moduli space $\P_\beta$ is the 
$(\sum_{i = 1}^n w_i)b + d$-dimensional  weighted projective
space: 
\[ \P(\underbrace{w_1, \ldots, w_1}_{b+1}, \ldots, 
\underbrace{w_n, \ldots, w_n}_{b+1}) . \]
It is easy to see that the Morrison-Plesser class defined by
the nef-partition is 
\[
  \Phi_\beta = (d_1[D_0])^{d_1 b}\cdots (d_r[D_0])^{d_r b}.
\]
Using 
$\<[D_0]^{\dim \P_\beta}\>_\beta = 1/w_1^{w_1 b + 1}\cdots w_n^{w_n b + 1}$, 
we obtain 

\begin{eqnarray*}
  I_P(y) &=& \sum_{b\ge 0} 
  \<P([D_1],\ldots,[D_n])(d_1[D_0])^{d_1 b}\cdots (d_r[D_0])^{d_r b}\>_\beta \,
  y^b 
  \\ &=& 
  P(w_1,\ldots,w_n)\sum_{b\ge 0} (d_1^{d_1} \cdots d_r^{d_r})^b
  \<[D_0]^{\dim\P_\beta}\>_\beta \, y^b 
  \\ &=&
  P(w_1,\ldots,w_n)\sum_{b\ge 0} (d_1^{d_1} \cdots d_r^{d_r})^b
  \frac{1}{w_1^{w_1 b + 1}\cdots w_n^{w_n b + 1}} \, y^b 
  \\ &=&
  \nu \cdot P(w_1,\ldots,w_n) \sum_{b\ge 0} \mu^b \, y^b
  \\ &=& 
  \frac{\nu \cdot P(w_1,\ldots,w_n)}{1 - \mu y}.
\end{eqnarray*}
\end{proof}

The convex hull of the vectors $v_1,\ldots,v_n$ is a reflexive polytope 
(simplex) $\nabla^* \subset M_\R\cong \R^d$. 
Let $\M :=  M \oplus \Z^r $ be an extension of the lattice $M$ and 
$\{ b_1,\ldots,b_r\}$  the standard basis of $\Z^r$.
The  $(d + r - 1)$-dimensional Cayley polytope 
\[
  \D = \Delta_1 * \cdots * \Delta_r 
\]
is the convex hull of 
$(d + r + 1)$ points:
$e_1 = (0,b_1),\ldots,e_r = (0,b_r)$ and 
$u_k = (v_k,b_j)$ $(k = 1,\ldots,d + 1)$, 
where  $v_k\in A_j$. 
 We denote this set of points by $\widetilde{A}$.
The points from $ \widetilde{A}$ are affinely dependent, 
while any proper subset 
of $ \widetilde{A}$ is affinely independent, i.e., defines a circuit
(see \cite[Chapter~7]{GKZ}). 
It is easy to see that the only affine relation (up to a real multiple) 
between the points from $\widetilde{A}$ is
\[
  d_1 e_1 + \cdots + d_r e_r -
  w_1 u_1 - \cdots - w_n u_n = 0.
\]
Thus by \cite[Chapter~7, Proposition~1.2]{GKZ}, polytope 
 $\D$ has exactly 
two triangulations: the triangulation ${\mathcal T} = {\mathcal T}_1$ 
with the simplices 
$conv(\widetilde{A} \setminus \{ e_i\})$, $i = 1,\ldots,r$, and  the 
triangulation ${\mathcal T}_2$ with the simplices 
$conv(\widetilde{A}  \setminus \{u_k\} )$, $k = 1,\ldots,n$. Note that
\begin{eqnarray}
\label{vol1}
  {\rm Vol}(conv(\widetilde{A} \setminus \{ e_i\})) = d_i,\quad i = 1,\ldots,r,
\end{eqnarray}
\begin{eqnarray}
\label{vol2}
  {\rm Vol}(conv(\widetilde{A}  \setminus \{ v_k\})) = w_k,\quad k = 1,\ldots,n.
\end{eqnarray}
Therefore ${\rm Vol}(\D) = \sum_{i=1}^r d_i = \sum_{k=1}^n w_k$.

Let 
\[
  f_j(t) := 1 - \sum_{i\, : \, v_i \in A_j} a_i t^{v_i}\in
  \C[t_1^{\pm 1},\ldots,t_d^{\pm 1}],\quad j = 1,\ldots,r
\]
be generic Laurent polynomials. Denote by
\[
  F(t) = t_{d + 1} f_1(t) + \cdots + t_{d + r} f_r(t)
\]
a Laurent polynomial whose support polytope is $\D$.

The next statement 
 follows directly from \cite[Chapter~9, Proposition~1.8]{GKZ} and 
from the equalities (\ref{vol1}), (\ref{vol2}).

\begin{prop}
The $\widetilde{A}$-discriminant of $F$ is equal (up to sign) to the binomial
\[
  D_{\widetilde{A}}(F) = 
  \prod_{k = 1}^n w_k^{w_k} - 
  \prod_{i = 1}^r d_i^{d_i}\prod_{k = 1}^n a_k^{w_k} = 
\prod_{k = 1}^n w_k^{w_k}(1 - \mu y) 
,
\]
where $y = \prod_{k = 1}^n a_k^{w_k}$ and 
the first summand in $D_{\widetilde{A}}(F)$ corresponds to the triangulation 
${\mathcal T}$.
\end{prop}

\begin{theo}
 Let $P(x_1,\ldots,x_n)\in \C[x_1,\ldots,x_n]$ be a
 $\overline{k}$-homogeneous polynomial  with  $|\overline{k}| = d$. 
Then 
\[
  R_P(a) =(-1)^d {\rm Res}_F^{\overline{k}} \left(
   t_{d+1}^{\overline{k_1}}\cdots t_{d+r}^{\overline{k_r}}
P(a_1 t^{v_1},\ldots, a_n t^{v_n})\right) = 
 \frac{\nu\cdot P(w_1,\ldots,w_n)}{1 - \mu y},
\]
where 
 $y := a_1^{w_1}\cdots a_n^{w_n}$. 
\end{theo}

\begin{proof} 
By Proposition~\ref{sum_res} the toric residue $R_P(a)$ is the following sum 
over the critical points $\xi$ of the polynomial $F_1(t,y):= 
f_1(t) + y_2 f_2(t) + \cdots + y_r f_r(t)$, where $(t,y)\in 
 (\C^*)^d \times (\C^*)^{r - 1}$:
\[
  R_P(a) = (-1)^d \sum_{\xi \in V_F}
  \frac{P(a_1\xi^{v_1}, \ldots,a_n\xi^{v_n} )}{F_1(\xi) H_{F_1}^1(\xi)}.  
\]
We rewrite the polynomial $F_1$ as 
\[
  F_1 = y_2 + \cdots + y_r + 1 - \sum_{i = 1}^n c_i t^{v_i},
\]
where $c_i = y_j\cdot a_i$ if $v_i\in A_j$.
Then at the critical point $\xi$, we have 
\[
  c_1 \frac{\xi^{v_1}}{w_1} = \cdots = c_n \frac{\xi^{v_n}}{w_n} = z
\]
and 
\[
  z^{w_1 + \cdots + w_n} = 
  \left(\frac{c_1}{w_1}\right)^{w_1} \cdots \left(\frac{c_n}{w_n}\right)^{w_n}.
\]
Moreover, at the critical points one has: 
\[
  f_2(\xi) = \cdots = f_r(\xi) = 0,
\]
which is equivalent to 
\[
  f_j(\xi) = 1 - \sum_{i:v_i\in A_j} a_i \xi^{v_i} = 
  1 - \left(\sum_{i:v_i\in A_j} w_i\right) \frac{z}{\eta_j} = 
  1 - \frac{d_j z}{\eta_j} = 0,\quad j = 2,\ldots,r,
\]
where $\eta_j$ is the value of $y_j$ at the critical point.
Hence, it is easy to see that $\eta_j = d_j z$, $(j = 2,\ldots,r)$ and 
$F = 1 - d_1 z$, which implies
\begin{equation}
\label{roots}
  z^{w_1 + \cdots + w_n} = 
  \left(\frac{a_1}{w_1}\right)^{w_1} \cdots 
  \left(\frac{a_n}{w_n}\right)^{w_n} d_2^{d_2}\cdots d_r^{d_r}\cdot
  z^{d_2 + \cdots + d_r},
\end{equation}
or, equivalently,
\[
  z^{d_1} = 
  \left(\frac{a_1}{w_1}\right)^{w_1} \cdots 
  \left(\frac{a_n}{w_n}\right)^{w_n} d_2^{d_2}\cdots d_r^{d_r}.
\]
The value of the Hessian $H_F^1$ at $\xi$
equals 
\[
  H_F^1(\xi) = (-1)^d w_1\cdots w_n d_1 \cdots d_r z^{d + r - 1}.
\]
Since there are exactly $w_1 + \cdots + w_n = d_1 + \cdots + d_r$
critical points of $F$, the summation over the critical points is equivalent 
to the summation over the roots of (\ref{roots}), we get
\begin{eqnarray*}
  && R_P(y)  =
  \sum_{z^{d_1} = 
  \left(\frac{a_1}{w_1}\right)^{w_1} \cdots 
  \left(\frac{a_n}{w_n}\right)^{w_n} d_2^{d_2}\cdots d_r^{d_r}}
  \frac{P(w_1,\ldots,w_n)}
  {w_1\cdots w_n d_1 (1 - d_1 z)} 
   \\ && =
  \frac{P(w_1,\ldots,w_n)}{w_1\cdots w_n} 
  \sum_{b \ge 0} 
  \left(\displaystyle d_1^{d_1}\cdots d_r^{d_r}\right)^b
  \left(\frac{a_1^{w_1}\cdots a_n^{w_n}}{w_1^{w_1}\cdots w_n^{w_n}}\right)^b
  \\ && = \frac{\nu\cdot P(w_1,\ldots,w_n)}{1 - \mu\, y}.
\end{eqnarray*}

\end{proof}

\section{Complete intersections in product of projective spaces}

In this section we check the Toric Residue Mirror Conjecture for 
nef-partitions corresponding to mirrors of complete intersections in 
a product of projective spaces
$\P = \P^{d_1}\times\cdots\times\P^{d_p}$ of dimension 
$d = d_1 + \cdots + d_p$. We set $n_i := d_i + 1$ and denote 
by $N = (n_{ij})$ 
an integral $p\times r$-matrix with non-negative elements having 
columns $\nn_1,\ldots,\nn_r\in\Z_{\ge 0}^p$. A complete intersection $V$ of
$r$ hypersurfaces $V_1,\ldots,V_r$ in $\P$ of multidegrees $\nn_1,\ldots,\nn_r$
is a Calabi-Yau $(d - r)$-fold if and only if 
$\sum_{j = 1}^r n_{ij} = n_i$ ($i = 1,\ldots,p$). 
We will use the standard  notation
\[
  \left(
  \begin{array}{c}
  \P^{d_1} \\
  \vdots   \\
  \P^{d_p} \\
  \end{array}
  \right|
  \left|
  \begin{array}{ccc}
  n_{11}  &\cdots & n_{1r} \\
   \vdots &       & \vdots\\
   n_{p1} &\cdots & n_{pr}
  \end{array}
  \right)
\]
to denote this complete intersection.

The cone of effective curves $K_{\rm eff}(\P)$ is isomorphic to 
$\R_{\ge 0}^p$ and its 
integral part $K_{\rm eff}(\P)_\Z$ consists of the points 
$\beta = (b_1,\ldots,b_p)\in\Z_{\ge 0}^p$. Thus, the Morrison-Plesser 
moduli spaces are the products of projective spaces: 
$\P_\beta = \P^{n_1 b_1 + d_1}\times\cdots\times\P^{n_p b_p + d_p}$ and 
the generating function for intersection numbers may be written
\[
  I_P(y) = \sum_{b_1,\ldots,b_p\ge 0} I(P,\beta)\,y_1^{b_1}\cdots y_p^{b_p}.
\]

\begin{theo} 
\label{I_RES}
The generating function for intersection numbers associated with the 
monomial $x^k = x_1^{k_1}\cdots x_p^{k_p}$ can be written as the integral
\begin{equation}
\label{I_as_Res}
  I_{x^k}(y) = \left(\frac{1}{2\pi i}\right)^p
  \int_\Gamma
  \frac{z_1^{k_1}\cdots z_p^{k_p} dz_1\wedge\cdots\wedge dz_p}
  {G_1(z)\cdots G_p(z)},
\end{equation}
where the polynomials $G_i$ have the form
\[
  G_i = z_i^{n_i} - 
  \prod_{j = 1}^r (n_{1j} z_1 + \cdots + n_{pj} z_p)^{n_{ij}} y_i,\quad
  i = 1,\ldots,p,
\]
and $\Gamma$ is the compact cycle in $\C^p$ defined by 
$\Gamma = \{|G_1| = \cdots = |G_p| = \delta \}$ for 
small positive $\delta$.
\end{theo}

\begin{proof} Let $[H_i]$ denote the class of a hyperplane section in 
$\P^{d_i}$. The class of the divisor $E_j$ defining hypersurface 
$V_j$ equals
\[
  [E_j] = n_{1j} [H_1] + \cdots + n_{pj} [H_p],\quad j  = 1,\ldots,r.
\]
Hence, the coefficients of the series $I_{x^k}(y)$ are 
\[
  \<[H_1]^{k_1}\cdots [H_p]^{k_p}
  \prod_{j = 1}^r 
  (n_{1j} [H_1] + \cdots + n_{pj} [H_p])^{n_{1j} b_1 + \cdots + n_{pj} b_p}
  \>_\beta.
\]
The lattice points $\beta = (b_1,\ldots,b_p)\in \Z_{\ge 0}^r$ in the 
integral part of the Mori cone $K_{\rm eff}(\P)$ correspond to the $p$ linear
relations  
\[
  b_i v_{i1} + \cdots + b_i v_{in_i} = 0,\quad i = 1,\ldots,p,
\]
where $v_{j1},\ldots,v_{jn_j}$ generate a lattice $M_j$ of rank $d_j$
($1\le j\le p$).
Therefore we set  $y_i := a_{i1}\cdots a_{in_i}$. 
Using the property of the integral 
\[
  \left(\frac{1}{2\pi i}\right)^p
  \int_{\gamma_\rho} 
  z_1^{m_1 - 1}\cdots z_p^{m_p - 1}\, dz = 
  \left\{
  \begin{array}{ll}
  1,   & m_1 = \cdots = m_p = 0, \\
  0,   & {\rm otherwise},
  \end{array}
  \right.
\]
where $\gamma_\rho = \{|z_1| = \cdots = |z_p| = \rho\}$
is the cycle winding around the origin ($\rho > 0$ is small) 
and the fact that the intersection numbers on $\P$ are 
\[
  \<[H_1]^{l_1}\cdots [H_p]^{l_p}\>_\beta = 
  \left\{
  \begin{array}{ll}
  1,   & l_j = n_j b_j + d_j, \quad j = 1,\ldots,r, \\
  0,   & {\rm otherwise},
  \end{array}
  \right.
\]
we can represent the functions $I(x^k,\beta)$ by integrals
\[
  I(x^k,\beta) = 
  \left(\frac{1}{2\pi i}\right)^p
  \int_{\gamma_\rho}
  \frac{z_1^{k_1}\cdots z_p^{k_p} \prod_{j = 1}^r
  (n_{1j} z_1 + \cdots + n_{pj} z_p)^{n_{1j} b_1 + \cdots + n_{pj} b_p}
  \, dz}{z_1^{n_1 b_1 + d_1 + 1}\cdots z_p^{n_p b_p + d_p + 1}}.
\]
Denote 
\[
  F_i(z) := \prod_{j = 1}^r
  (n_{1j} z_1 + \cdots + n_{pj} z_p)^{n_{ij}},\quad i = 1,\ldots,p.
\]
If $z\in\gamma_\rho$ for some fixed $\rho$, then the geometric series
\[
  z_1^{k_1 - n_1}\cdots z_p^{k_p - n_p}
  \sum_{b_1,\ldots,b_p\ge 0}
  \left(\frac{F_1(z) y_1}{z_1^{n_1}} \right)^{b_1}\cdots
  \left(\frac{F_p(z) y_p}{z_p^{n_p}} \right)^{b_p} = 
  \frac{z_1^{k_1}\cdots z_p^{k_p}}
  {\prod_{i = 1}^p \left[z_i^{n_i} - F_i(z) y_i\right]}
\]
converges absolutely and uniformly for all $y$ from the neighborhood 
${\mathcal U}_\varepsilon = \{y:||y|| < \varepsilon\}$, where 
$0 < \varepsilon < \min_{i = 1,\ldots,p}(\rho^{n_i}/M_i)$,
$M_i = \max_{z\in\gamma_\rho} |F_i(z)|$.
Integrating the last expression and changing the order of integration 
and summation, we get
\[
  I_{x^k}(y) = 
  \left(\frac{1}{2\pi i}\right)^p \int_{\gamma_\rho}
  \frac{z_1^{k_1}\cdots z_p^{k_p}\,dz_1\wedge\cdots\wedge dz_p}
  {\prod_{i = 1}^p \left[z_i^{n_i} - F_i(z) y_i\right]}.
\]
For fixed $y\in {\mathcal U}_\varepsilon$, using the  
Rouch\'e's principle for residues 
(see \cite[Chapter~2, \S 8]{Tsikh} or \cite[Lemma~4.9]{Aizenberg-Yuzhakov}), 
the cycle $\gamma_\rho$ can be replaced by another homologically 
equivalent cycle $\Gamma$ of the following form: 
\[
  \gamma_\rho\sim\Gamma = 
  \{z:|z_1^{n_1} - F_1(z) y_1| = \cdots = |z_p^{n_p} - F_p(z) y_p| = \delta\}
\]
for some $\delta > 0$. 
Therefore, we have
\[
  I_{x^k}(y) = 
  \left(\frac{1}{2\pi i}\right)^p \int_\Gamma
  \frac{z_1^{k_1}\cdots z_p^{k_p}\,dz_1\wedge\cdots\wedge dz_p}
  {\prod_{i = 1}^p \left[z_i^{n_i} - F_i(z) y_i\right]}
\]
which finishes the proof.
\end{proof}

The Conjecture 4.6 follows now from a general result in 
\cite{Batyrev-Materov2} which identifies 
\[  \left(\frac{1}{2\pi i}\right)^p
  \int_\Gamma
  \frac{z_1^{k_1}\cdots z_p^{k_p} dz_1\wedge\cdots\wedge dz_p}
  {G_1(z)\cdots G_p(z)} \] 
with the toric residue.

\section{Computation of Yukawa $(d - r)$-point functions}

Let $\Delta = \Delta_1 + \cdots + \Delta_r$ be 
a nef-partition of a reflexive polytope $\Delta$, 
$A_i \subset  \partial \nabla^* \cap \Delta_i \cap M$ a subset 
containing all nonzero vertices of  $\Delta_i$ $( 1\leq i \leq r)$.
We set $A_1 \cup \cdots \cup A_r:= \{ v_1, \ldots, v_n \}$ and   
consider a $\Delta_1 * \cdots * \Delta_r$-regular sequence of Laurent
polynomials
\[
  f_j(t) := 1 - \sum_{i: v_i\in A_j} a_i t^{v_i}
  \in \C[t_1^{\pm 1},\ldots,t_d^{\pm 1}],\quad j = 1,\ldots,r, 
\]
which define $r$ affine hypersurfaces
\[ Z_{f_j} := \{ t \in \T \cong (\C^*)^d \; : \; f_j(t) = 0\}, 
\quad j = 1,\ldots,r, \]
The compactification $\overline{Z}_f$ in $\P_{\Delta}$ 
of the affine  complete intersection  $Z_f := Z_{f_1}\cap\cdots\cap Z_{f_r}$
is a $(d-r)$-dimensional 
projective Calabi-Yau variety with at worst Gorenstein canonical 
singularities. Using the Poincar\'e residue mapping
\[
  {\bf Res}\,:\, H^d(\T\setminus Z_{f_1}\cup\cdots\cup Z_{f_r})\rightarrow
  H^{d - r}(Z_{f_1}\cap\cdots\cap Z_{f_r})
\]
one can construct a nowhere vanishing section of the canonical bundle 
of $\overline{Z}_f$ as 
\[
  \Omega := {\bf Res}
  \left(
  \frac{1}{f_1\cdots f_r}
  \frac{d t_1}{t_1}\wedge\cdots\wedge\frac{d t_d}{t_d}
  \right).
\]

\begin{dfn} 
Let $Q(x_1,\ldots,x_n)\in\Q[x_1,\ldots,x_n]$ be a homogeneous polynomial 
of degree $d - r$. 
The $Q$-{\it Yukawa  $(d-r)$-point function} is defined by the formula 
\begin{equation*}
\label{def_of_Yuk}
  Y_Q(a_1,\ldots,a_n) := 
  \frac{(-1)^{\frac{(d - r)(d - r - 1)}{2}}}{(2\pi i)^{d - r}}
  \int_{Z_f} \Omega \wedge 
  \displaystyle
  Q\left(a_1 \frac{\partial}{\partial a_1},\ldots,
         a_n \frac{\partial}{\partial a_n}\right) \Omega,
\end{equation*}
where the differential operators 
$a_1\partial/\partial a_1,\ldots,a_n\partial/\partial a_n$ are determined  by
the Gau\ss-Manin connection.
If $\overline{k} = (\overline{k_1},\ldots, \overline{k_r})$ is a nonnegative
integral vector with $|\overline{k}| = d-r$ and $Q(x_1,\ldots,x_n)$ is a 
$\overline{k}$-homogeneous polynomial (${\rm deg}\, x_j = k_i$ 
$\Leftrightarrow$ $v_j \in A_i$), then 
\[  Q\left(a_1 \frac{\partial}{\partial a_1},\ldots,
         a_n \frac{\partial}{\partial a_n}\right) \Omega = (-1)^{d - r}
{\bf Res}
  \left(
  \frac{Q(a_1t^{v_1}, \ldots,a_nt^{v_n}) }{f_1^{\overline{k_1}+1}\cdots 
f_r^{\overline{k_r}+1}}
  \frac{d t_1}{t_1}\wedge\cdots\wedge\frac{d t_d}{t_d}
  \right).
\] 
\end{dfn}

\begin{theo} Let $Q(x_1,\ldots,x_n)\in \C[x_1,\ldots,x_n]$ be a
 $\overline{k}$-homogeneous polynomial  with  $|\overline{k}| = d-r$.
 We define 
\[
  P(x_1,\ldots,x_n) := 
  \prod_{j = 1}^r \left(\sum_{i:v_i\in A_j} x_i\right)
  Q(x_1,\ldots,x_n).
\]
Then the Yukawa $(d - r)$-point function is equal to the $k$-mixed 
toric residue
\[
  Y_Q(a_1,\ldots,a_n) = 
  (-1)^d {\rm Res}_F^{k} \left(
   t_{d+1}^{\overline{k_1}+1}\cdots t_{d+r}^{\overline{k_r}+1}
P(a_1 t^{v_1},\ldots, a_n t^{v_n})\right).
\]
\end{theo}

\begin{proof} We sketch only the idea of the proof.  
The hypersurface $$Z_F = \{
t_{d+1} f_1 + \cdots + t_{d+r}f_r = 1\}$$ in $(\C^*)^d \times 
\C^r$ is a  
$\C^{r-1}$-bundle over 
$(\C^*)^d \setminus \left(Z_{f_1} \cap \cdots \cap Z_{f_r} \right)$. 
This fact allows one to identify the primitive parts of the cohomology groups 
$H^{d-r}( Z_{f_1} \cap \cdots \cap Z_{f_r})$ and $H^{d-1}(Z_F)$ together 
with their intersection forms. By the result of Mavlyutov 
\cite{Mavlyutov2}, one can compute the intersection form on $H^{d-1}(Z_F)$
using toric resides.   
\end{proof} 

\begin{exam} 
Consider the mirror family $V^*$ to Calabi-Yau  complete intersections   $V$ 
of $r$ hypersurfaces of degrees $d_1,\ldots,d_r$ respectively in $\P^d$, 
$d_1 + \cdots + d_r = d+1$. Its  nef-partition can be constructed 
as follows.  
Let $v_1,\ldots,v_d$ be a basis vectors of the lattice $M$ and 
\[
  v_{d + 1} := -v_1 - \cdots - v_d.
\]
We divide the set $\{v_1,\ldots,v_{d + 1}\}$ into a disjoint union of $r$
subsets $A_1,\ldots,A_r$ such that $|A_i| = d_i$. 
 For $j = 1,\ldots,r$, we define Laurent polynomials 
\[
  f_j(t) := 1 - \sum_{i : v_i\in A_j} a_i t^{v_i}
  \in \C[t_1^{\pm 1},\ldots,t_d^{\pm 1}].
\]
Then the affine part of $V^*$ is the complete intersection $Z_f\subset\T$ 
of hypersurfaces 
$Z_{f_1},\ldots,Z_{f_r}\subset\T$ defined by polynomials $f_1,\ldots,f_r$.
The Yukawa coupling for $V^*$  has been computed in  
\cite[Proposition~5.1.2]{BvS}: 
\[
  Y_Q(y) = \frac{d_1\cdots d_r Q(1,\ldots,1)}{1 - \mu y},
\]
where $y = a_1\cdots a_n$ and $\mu = \prod_{i = 1}^r d_i^{d_i}$.
\end{exam}

\begin{exam} 
Consider  Calabi-Yau varieties $V$ obtained as the complete intersection 
of hypersurfaces $V_1,V_2,V_3$
in $\P^3\times\P^3$ of degrees $(3,0)$, $(0,3)$ and $(1,1)$ respectively
of the type 
\[
 \left(
 \begin{array}{c}
 \P^3 \\
 \P^3 
 \end{array}
 \right|
 \left|
 \begin{array}{ccc}
 3 & 0 & 1 \\
 0 & 3 & 1 
 \end{array}
 \right).
\]
Let $M\cong\Z^6$ and $\nabla^* = 
conv(\Delta_1\cup\Delta_2\cup\Delta_3)\subset M_\R$
be a reflexive polytope  defined by the polytopes
$\Delta_1 := conv\{0,v_1,v_2,v_3\}$,
$\Delta_2 := conv\{0,v_5,v_6,v_7\}$ and 
$\Delta_3 := conv\{0,v_4,v_8\}$, where
\begin{eqnarray*} 
  &&v_1 = (1,0,0,0,0,0),\,v_2 = (0,1,0,0,0,0),\,v_3 = (0,0,1,0,0,0),\\
  &&v_4 = (-1,-1,-1,0,0,0),\,v_5 = (0,0,0,1,0,0),\,v_6 = (0,0,0,0,1,0),\\
  &&v_7 = (0,0,0,0,0,1),\,v_8 = (0,0,0,-1,-1,-1).
\end{eqnarray*}
The  nef-partition 
$\Delta = \Delta_1 + \Delta_2 + \Delta_3$ corresponds to mirrors $V^*$ 
of  $V = V_1 \cap V_2 \cap V_3$.
We define  the disjoint sets: $A_1 := \{v_1,v_2,v_3\}$, 
$A_2 := \{v_5,v_6,v_7\}$, $A_3 := \{v_4,v_8\}$ and  the Laurent
polynomials 
\begin{eqnarray*} 
  &&f_1(t) := 1 - \sum_{i:v_i\in A_1} a_i t^{v_i} = 
  1 - a_1 t_1 - a_2 t_2 - a_3 t_3, \\
  &&f_2(t) := 1 - \sum_{i:v_i\in A_2} a_i t^{v_i} = 
  1 - a_5 t_4 - a_6 t_5 - a_7 t_6, \\
  &&f_3(t) := 1 - \sum_{i:v_i\in A_3} a_i t^{v_i} = 
  1 - a_4 t_1^{-1} t_2^{-1} t_3^{-1} - a_8 t_4^{-1} t_5^{-1} t_6^{-1}.
\end{eqnarray*}
The complete intersection $Z_f := Z_{f_1}\cap Z_{f_2}\cap Z_{f_3}$ of the
affine hypersurfaces \[ Z_{f_j} = \{t\in (\C^*)^6\,:\,f_j(t) = 0\},
\quad j = 1,2,3 \] is an affine part of 
$V^*$. 

Denote by $y_1 = 3^3 a_1 a_2 a_3 a_4$, $y_2 =  3^3 a_5 a_6 a_7 a_8$ the
new variables and by $\theta_1 := y_1\partial/\partial y_1$,
$\theta_2 := y_2\partial/\partial y_2$ the corresponding logarithmic partial
derivations. Given a form-residue
\[
  \Omega := {\bf Res}\left(\frac{1}{f_1 f_2 f_3}\frac{d t_1}{t_1}\wedge
  \cdots\wedge\frac{d t_6}{t_6}\right)\in H^3(Z_f),
\]
we compute the $2$-parameter Yukawa couplings defined as integrals
\[
  Y^{(k_1,k_2)}(y_1,y_2) = \frac{-1}{(2\pi i)^3 }
  \int_{Z_f} \Omega\wedge\theta_1^{k_1}\theta_2^{k_2}\Omega,\quad
  k_1 + k_2 = 3.
\]

\begin{prop} 
\label{BvS_Yuk}
The Yukawa couplings are
\begin{eqnarray*}
  Y^{(3,0)}(y_1,y_2) = \frac{9 y_1}{(1 - y_1 - y_2)(1 - y_1)^2},\;
  Y^{(2,1)}(y_1,y_2) = \frac{9}{(1 - y_1 - y_2)(1 - y_1)},\\
  Y^{(1,2)}(y_1,y_2) = \frac{9}{(1 - y_1 - y_2)(1 - y_2)},\;
  Y^{(0,3)}(y_1,y_2) = \frac{9 y_2}{(1 - y_1 - y_2)(1 - y_2)^2}.
\end{eqnarray*}
\end{prop}

\begin{rem}
Note that the functions in Proposition~\ref{BvS_Yuk} are completely consistent
with Yukawa couplings from \cite[\S 8.3]{BvS}. Indeed, if we put 
$K(y_1,y_2) = Y^{(3,0)} + 3 Y^{(2,1)} + 3 Y^{(1,2)} + Y^{(0,3)}$ 
and consider restriction
to the diagonal subfamily $y = y_1 = y_2$, 
then we get the same expression as in \cite{BvS}:
\[
  K(y,y) = \frac{18(3 - 2y)}{(1 - 2y)(1 - y)^2}.
\]
\end{rem}

Denote by $F(t) := t_7 f_1(t) + t_8 f_2(t) + t_9 f_3(t)$ the Cayley polynomial 
associated with Laurent polynomials $f_1, f_2, f_3$, and by 
$\D = \Delta_1 * \Delta_2 * \Delta_3\subset \M_\R = M_\R\oplus \R^3$ 
its supporting polytope which is the
Cayley polytope associated with $\Delta_1,\Delta_2,\Delta_3$.

\begin{prop} Let $A := \D \cap \M$ and $F(t)$ be the Cayley polynomial as above.
Then the principal $A$-determinant of $F$ has the form
\[
  E_A(F) = (a_1\cdots a_8)^{12} (1 - y_1)^3(1 - y_1)^3(1 - y_1 - y_2).
\]
\end{prop}

\begin{rem} It is easy to see that the products of $Y^{(k_1,k_2)}$ by 
$E_A(F)$ are polynomials in $a_1,\ldots,a_8$.
\end{rem}

Let us find the generating function $I_P(y)$ for the monomial 
$P(x) = x_1^{k_1} x_2^{k_2}$. There are two linear independent integral
relations between $v_1,\ldots,v_8$:
\[
  v_1 + \cdots + v_4 = 0,\quad v_5 + \cdots + v_8 = 0.
\]
Hence the Mori cone $K_{\rm eff}(\P)$ is spanned by the vectors 
\[
  l^{(1)} = (1,1,1,1,0,0,0,0),\quad l^{(2)} = (0,0,0,0,1,1,1,1)
\]
and the Morrison-Plesser moduli spaces are 
$\P_\beta = \P^{4 b_1 + 3}\times\P^{4 b_2 + 3}$ 
($b_1, b_2\in \Z_{\ge 0}$). 
The cohomology ring of $\P_\beta$ is
generated by two hyperplane classes: $[H_1]$ and $[H_2]$.
We set $E_1 := 3[H_1]$, $ E_2 := 3[H_2]$ and $E_3 := [H_1] + [H_2]$.
Then the nef-partition of the 
anticanonical divisor $-K_\P$ corresponding to the nef-partition 
$\Delta = \Delta_1 + \Delta_2 + \Delta_3$ is defined by 
$-K_\P = E_1 + E_2 + E_3$. Therefore, the Morrison-Plesser cohomology class 
associated with the nef-partition of $-K_\P$ equals
\[
  \Phi_\beta = (3[H_1])^{3 b_1} (3[H_2])^{3 b_2} 
  ([H_1] + [H_2])^{b_1 + b_2}
\]
and the generating function for intersection numbers can be written
\[
  I_P(y) = \sum_{b_1,b_2 \ge 0}
  \<[H_1]^{k_1 + 3 b_1 + 1}[H_2]^{k_2 + 3 b_2 + 1}
  ([H_1] + [H_2])^{b_1 + b_2 + 1}\>_\beta \, y_1^{b_1} y_2^{b_2}.
\]
Intersection theory on $\P_\beta$ implies
\[
  I_P(y) = 9 \sum_{b_1,b_2 \ge 0} 
  \frac{(b_1 + b_2 + 1)!}{(b_1 - k_1 + 2)!(b_2 - k_2 + 2)!}\,
  y_1^{b_1} y_2^{b_2}.
\]
By Theorem~\ref{I_RES} we can write  $I_P(y)$ as the 
integral
\[
  I_P(y) = \frac{1}{(2\pi i)^2}\int_\Gamma
  \frac{9 z_1^{k_1 + 1} z_2^{k_2 + 1}(z_1 + z_2)\,d z_1\wedge d z_2}
  {(z_1^4 - z_1^3(z_1 + z_2)y_1)(z_2^4 - z_2^3(z_1 + z_2)y_2)}
\]
with the cycle $\Gamma = \{(z_1,z_2)\in\C^2\,:\,
|z_1^4 - z_1^3(z_1 + z_2)y_1| = 
|z_2^4 - z_2^3(z_1 + z_2)y_2| = \delta\}$, 
$\delta > 0$. Computing the last integrals, we get 
the same rational functions as in Proposition~\ref{BvS_Yuk}.
\end{exam}

\begin{exam}
Consider  Calabi-Yau varieties  $V$ obtained as  complete
intersections of two hypersurfaces of degrees $(4,0)$, $(1,2)$ in 
$\P = \P^4\times\P^1$ which corresponds to the configuration 
\[
 \left(
 \begin{array}{c}
 \P^4 \\
 \P^1 
 \end{array}
 \right|
 \left|
 \begin{array}{cc}
 4 & 1  \\
 0 & 2  
 \end{array}
 \right).
\]
This example was investigated in details by Hosono, Klemm, Theisen 
and Yau (cf. \cite{Yau}).
The corresponding nef-partition $\Delta = \Delta_1 + \Delta_2
\subset M_\R\cong \R^5$ consists of polytopes $\Delta_1 := 
conv\{0, v_1,v_2,v_3,v_4\}$ 
and $\Delta_2 := conv\{0, v_5,v_6,v_7\}$, where 
\begin{eqnarray*}
  &&v_1 = (1,0,0,0,0),\quad v_2 = (0,1,0,0,0),\quad v_3 = (0,0,1,0,0),\quad
    v_4 = (0,0,0,1,0),\quad \\
  &&v_5 = (-1,-1,-1,-1,0),\quad v_6 = (0,0,0,0,1),\quad
    v_7 = (0,0,0,0,-1).
\end{eqnarray*}
We have two  disjoint sets: 
$A_1 := \{v_1,v_2,v_3,v_4\}$ and $A_2 := \{v_5,v_6,v_7\}$ 
which are the vertices of the reflexive polytope $\nabla^*$ and
define
the Laurent polynomials 
\begin{eqnarray*}
  &&f_1(t) = 1 - \sum_{i:v_i\in A_1} a_i t^{v_i} = 
  1 - a_1 t_1 - a_2 t_2 - a_3 t_3 - a_4 t_4,\\
  &&f_2(t) = 1 - \sum_{i:v_i\in A_2} a_i t^{v_i} = 
  1 - a_5 (t_1 t_2 t_3 t_4)^{-1} - a_6 t_5 - a_7 t_5^{-1}.
\end{eqnarray*}

Denote by 
$y_1 := a_1 \cdots a_5$, 
$y_2 := a_6 a_7$ the new variables and 
$\theta_1 := y_1\partial/\partial y_1$,
$\theta_2 := y_2\partial/\partial y_2$ 
the corresponding logarithmic partial derivations.
Let $\Omega$ be a form defined by 
\[
  \Omega := {\bf Res}\left(\frac{1}{f_1 f_2}
  \frac{d t_1}{t_1}\wedge\frac{d t_2}{t_2}\right)\in 
  H^3(Z_{f_1}\cap Z_{f_2}).
\]
Then the Yukawa coupling associated with $f_1,f_2$ is the integral 
\[
  Y^{(k_1,k_2)}(y_1,y_2) := 
  \frac{-1}{(2\pi i)^3}\int_{Z_f}
  \Omega\wedge\theta_1^{k_1}\theta_2^{k_2}\Omega,
  \quad k_1 + k_2 = 3,
\]
where $Z_f := Z_{f_1}\cap Z_{f_2}$ is an affine Calabi-Yau complete
intersection which compactification forms a mirror dual family to $V$.

\begin{prop} \cite{Yau} 
\label{ex_Yuk}
The Yukawa couplings $Y^{(k_1,k_2)}(y)$ are:
\begin{eqnarray*}
  &&Y^{(3,0)}(y) = \frac{8}{D_0},\quad
    Y^{(2,1)}(y) = \frac{4(1 - 256 y_1 + 4 y_2)}{D_0 D_1},\\
  &&Y^{(1,2)}(y) = \frac{8 y_2(3 - 512 y_1 + 4 y_2)}{D_0 D_1^2},\\
  &&Y^{(0,3)}(y) = \frac{4 y_2(1 - 256 y_1 + 24 y_2 - 3072 y_1 y_2 + 16 y_2^2)}
                     {D_0 D_1^3},
\end{eqnarray*}
where
\[
  D_0 := (1 - 256 y_1)^2 - 4 y_2,\quad D_1 := 1 - 4 y_2.
\]
\end{prop}

\medskip 

Let $F(t) := t_6 f_1(t) + t_7 f_2(t)$ be the Cayley polynomial associated with 
$f_1(t)$ and $f_2(t)$. Its support polytope is the Cayley polytope 
$\D = \Delta_1*\Delta_2\subset\M_\R = M_\R \oplus \R^2$ which is the convex
hull of the vectors: 
\begin{eqnarray*}
  &&u_1 = (0,0,0,0,0;1,0),\,
  u_2 = (1,0,0,0,0;1,0),\,
  u_3 = (0,1,0,0,0;1,0),\\
  &&u_4 = (0,0,1,0,0;1,0),\,
  u_5 = (0,0,0,1,0;1,0),\,
  u_6 = (0,0,0,0,0;0,1),\\
  &&u_7 = (-1,-1,-1,-1,0;0,1),\,
  u_8 = (0,0,0,0,1;0,1),\,
  u_9 = (0,0,0,0,-1;0,1).
\end{eqnarray*}

\begin{prop} Let
$A := \{u_1,\ldots,u_9\}\subset \M$. Then the principal $A$-determinant of
$F(t)$  has the  following form: 
\begin{eqnarray*}
  &&E_A(F) = a_1^8 a_2^8 a_3^8 a_4^8 a_5^8 a_6^5 a_7^5
  D_0 D_1^4 = \\
  &&-640 a1^8 a2^8 a3^8 a4^8 a5^8 a6^8 a7^8 - 
  16777216 a1^{10} a2^{10} a3^{10} a4^{10} a5^{10} a6^8 a7^8 + \\
  &&160 a1^8 a2^8 a3^8 a4^8 a5^8 a6^7 a7^7 + 
  \underline{16777216 a1^{10} a2^{10} a3^{10} a4^{10} a5^{10} a6^9 a7^9} + \\
  &&\underline{65536 a1^{10} a2^{10} a3^{10} a4^{10} a5^{10} a6^5 a7^5} + 
  6291456 a1^{10} a2^{10} a3^{10} a4^{10} a5^{10} a6^7 a7^7 - \\
  &&1048576 a1^{10} a2^{10} a3^{10} a4^{10} a5^{10} a6^6 a7^6 - 
  \underline{1024 a1^8 a2^8 a3^8 a4^8 a5^8 a6^{10} a7^{10}} + \\
  &&1280 a1^8 a2^8 a3^8 a4^8 a5^8 a6^9 a7^9 - 
  512 a1^9 a2^9 a3^9 a4^9 a5^9 a6^5 a7^5 - \\
  &&20 a1^8 a2^8 a3^8 a4^8 a5^8 a6^6 a7^6 + 
  131072 a1^9 a2^9 a3^9 a4^9 a5^9 a6^8 a7^8 + \\
  &&\underline{a1^8 a2^8 a3^8 a4^8 a5^8 a6^5 a7^5} + 
  8192 a1^9 a2^9 a3^9 a4^9 a5^9 a6^6 a7^6 - \\
  &&49152 a1^9 a2^9 a3^9 a4^9 a5^9 a6^7 a7^7 - 
  131072 a1^9 a2^9 a3^9 a4^9 a5^9 a6^9 a7^9,
\end{eqnarray*}
where the terms corresponding to the vertices of the 
Newton polytope of $E_A(F)$ are
underlined. 
\end{prop}

\begin{proof}
The principal $A$-determinant can be found by using the algorithm 
proposed by A.~Dickenstein and B.~Sturmfels \cite{Dickenstein-Sturmfels} via
the computation of the corresponding Chow forms.
\end{proof}

\begin{rem} We note  that $D_0$ is the principal component of the 
discriminant locus 
$E_A(F) = 0$ and the component $D_1$ corresponds to the edge 
$\Gamma$ of $\D$ with 
\[
 \Gamma\cap\M = \{u_6,u_8,u_9\} = 
 \{(0,0,0,0,0;0,1),(0,0,0,0,1;0,1),(0,0,0,0,-1;0,1)\}.
\]
\end{rem} 

The Newton polytope of $E_A(F)$ is the secondary polytope ${\rm Sec}(A)$
depicted in Figure~\ref{Sec_P4_P1}. The vertices of ${\rm Sec}(A)$
are in one-to-one 
correspondence with coherent triangulations 
${\mathcal T}_1,\ldots,{\mathcal T}_4$ of $\D$ which are:

\begin{eqnarray*}
  {\mathcal T}_1 &=&
\{\<u_1,u_3,u_4,u_5,u_6,u_7,u_8\>,
  \<u_1,u_2,u_4,u_5,u_6,u_7,u_8\>,
  \<u_1,u_2,u_3,u_5,u_6,u_7,u_8\>,\\
&&\<u_1,u_2,u_3,u_4,u_6,u_7,u_8\>,
  \<u_1,u_2,u_3,u_4,u_5,u_6,u_8\>,
  \<u_1,u_3,u_4,u_5,u_6,u_7,u_9\>,\\
&&\<u_1,u_2,u_4,u_5,u_6,u_7,u_9\>,
  \<u_1,u_2,u_3,u_5,u_6,u_7,u_9\>, 
  \<u_1,u_2,u_3,u_4,u_6,u_7,u_9\>,\\
&&\<u_1,u_2,u_3,u_4,u_5,u_6,u_9\>\}.
\\
  {\mathcal T}_2 &=&
\{\<u_1,u_3,u_4,u_5,u_7,u_8,u_9\>,
  \<u_1,u_2,u_4,u_5,u_7,u_8,u_9\>,
  \<u_1,u_2,u_3,u_5,u_7,u_8,u_9\>, \\
&&\<u_1,u_2,u_3,u_4,u_7,u_8,u_9\>,
  \<u_1,u_2,u_3,u_4,u_5,u_8,u_9\>\}.
\\
  {\mathcal T}_3 &=&
\{\<u_2,u_3,u_4,u_5,u_7,u_8,u_9\>,
  \<u_1,u_2,u_3,u_4,u_5,u_7,u_9\>,
  \<u_1,u_2,u_3,u_4,u_5,u_7,u_8\>\}.
\\
  {\mathcal T}_4 &=&
\{\<u_1,u_2,u_3,u_4,u_5,u_7,u_8\>,
  \<u_1,u_2,u_3,u_4,u_5,u_7,u_9\>,
  \<u_2,u_3,u_4,u_5,u_6,u_7,u_8\>,\\
&&\<u_2,u_3,u_4,u_5,u_6,u_7,u_9\>\}.
\end{eqnarray*}

\begin{figure}[ht]
  \centering \includegraphics*[scale=0.5]{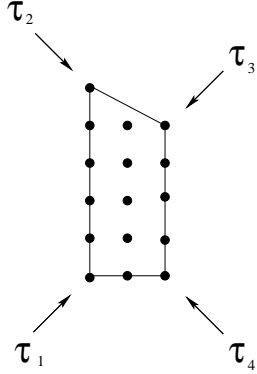}
  \caption{Secondary polytope and coherent triangulations}
  \label{Sec_P4_P1}
\end{figure}

\medskip

The generating function $I_P(y)$ for intersection numbers
corresponding to the monomial $P(x) = x_1^{k_1}x_2^{k_2}$
can be computed from the intersection theory on the Morrison-Plesser 
moduli spaces. Using  the two independent integral relations 
between $v_1,\ldots,v_7$
\[
  v_1 + \cdots + v_5 = 0,\quad v_6 + v_7 = 0,
\]
we see that the Mori cone $K_{\rm eff}(\P)$ is 
spanned by two vectors
\[
  l^{(1)} = (1,1,1,1,1,0,0),\quad 
  l^{(2)} = (0,0,0,0,0,1,1).
\]
The Morrison-Plesser 
moduli spaces are 
$\P_\beta = \P^{5b_1 + 4}\times \P^{2b_2 + 1}$ $(b_1, b_2 \in \Z_{\geq 0})$. 
The cohomology of 
$\P_\beta$ are generated by the hyperplane 
classes $[H_1]$ and $[H_2]$. Let $E_1 := [H_1] + 2[H_2]$ and $E_2 := 4[H_1]$. 
Then the nef-partition $\Delta = \Delta_1 + \Delta_2$ of polytopes induces the 
nef-partition of the anticanonical divisor 
\[
  -K_\P = E_1 + E_2 = ([H_1] + 2[H_2]) + (4[H_1]).
\]
It is straightforward to see that the corresponding Morrison-Plesser class is 
\[
  \Phi_\beta = ([H_1] + 2[H_2])^{b_1 + 2b_2}(4[H_1])^{4b_1}.
\]
So we get 
\[
  I_P(y) = \sum_{b_1,b_2\ge 0} 
  \< [H_1]^{k_1} [H_2]^{k_2} 
  ([H_1] + 2[H_2])^{b_1 + 2b_2 + 1}(4[H_1])^{4b_1 + 1}
  \>_\beta\, y_1^{b_1} y_2^{b_2}.
\]
Using the intersection theory on $\P_\beta$,  we obtain

\[
  I_P(y) = \sum_{b_1,b_2\ge 0} 
  2^{8 b_1 + 2 b_2 - k_2 + 3} 
  \frac{(b_1 + 2b_2 + 1)!}{(b_1 - k_1 + 3)!(2b_2 - k_2 + 1)!}\, 
  y_1^{b_1} y_2^{b_2}.
\]
By Theorem~\ref{I_RES} the function $I_P(y)$ admits the integral 
representation:
\[
  I_P(y) = \frac{1}{(2\pi i)^2} \int_\Gamma
  \frac{4 z_1^{k_1 + 1} z_2^{k_2}(z_1 + 2z_2)d z_1\wedge d z_2}
  {(z_1^5 - (z_1 + 2 z_2)(4 z_1)^4 y_1)
   (z_2^2 - (z_1 + 2 z_2)^2 y_2)}
\]
with the cycle \[\Gamma = \{(z_1,z_2)\in \C^2\,:\,
|z_1^5 - (z_1 + 2 z_2)(4 z_1)^4 y_1| = 
|z_2^2 - (z_1 + 2 z_2)^2 y_2| = \delta\},\]
where $\delta$ is positive.
These integrals can be easily computed and yield
the same rational functions as in Proposition~\ref{ex_Yuk}.
\end{exam}


\end{document}